\documentclass[11pt,leqno]{article}

\usepackage{amsmath,amsthm}
\usepackage {latexsym}
\usepackage{amssymb}

\newcommand{\bse}{\begin{equation}}

\newcommand{\bear}{\begin{eqnarray}}
\newcommand{\eear}{\end{eqnarray}}

\newcommand{\supp}{\mbox{\rm supp}}

\newcommand{\half}{\frac{1}{2}}

\newcommand{\eps}{{\varepsilon}}
\newcommand{\R}{{\mathbb R}}

\newcommand{\Compl}{{\mathbb C}}

\newcommand{\meas}{{\mathcal M}}

\newcommand{\sign}{\mbox{sign}}

\newcommand{\les}{\lesssim}

\newcommand{\Laplace}{\triangle}

\newcommand{\kato}{{\mathcal K}}

\newcommand{\la}{\lambda}

\def\Lap{\Delta}
\def\nn{\nonumber}

\newtheorem{theorem}{Theorem}
\newtheorem{lemma}[theorem]{Lemma}
\newtheorem{defi}[theorem]{Definition}
\newtheorem{cor}[theorem]{Corollary}

\newtheorem{proposition}[theorem]{Proposition}

\theoremstyle{remark}

\def\Pac{P_{a.c.}}
\def\la{\langle}
\def\ra{\rangle}

\def\norm[#1][#2]{\|#1\|_{#2}}

\def\y{{\bf y}}
\def\x{{\bf x}}
\def\Rnp{R_0^+(\lambda^2)}
\def\Bp{B^+}
\def\Btp{\tilde{B}^+}

\def\norm[#1][#2]{\Vert #1 \Vert_{#2}}
\def\impl{\Rightarrow}

\begin{document}

\title{Dispersive estimates for Schr\"odinger operators in
dimensions one and three}
\date{}

\author{M.\ Goldberg, W.\ Schlag}

\maketitle

\section{Introduction}

This paper deals with dispersive, i.e., $L^1(\R^d)\to L^\infty(\R^d)$
estimates for the time evolutions $e^{itH}P_{ac}$ where $H=-\Laplace+V$
and $P_{ac}$ is the projection onto the absolutely continuous spectral subspace. We 
restrict ourselves to the cases $d=1$ and $d=3$. Our goal is to assume
as little as possible on the potential $V=V(x)$ in terms of decay
or regularity. More precisely, we prove the following theorems.

\begin{theorem} 
\label{thm:thm1}
Let $V\in L^1_1(\R)$, i.e., $\int_{-\infty}^\infty 
|V(x)|(1+|x|)\,dx<\infty$,  and assume that there is no
resonance at zero energy. Then for all $t$
\begin{equation}
\label{eq:d1dec}  \big\|e^{itH} P_{ac}(H)\big\|_{1\to\infty} \les |t|^{-\half} 
\end{equation}
where $H=-\frac{d^2}{dx^2}+V$. 
The conclusion holds for all $V \in L^1_2(\R)$, i.e., $\int_{-\infty}^\infty 
|V(x)|(1+|x|)^{2}\,dx<\infty$,
whether or not there is a resonance at zero energy.
\end{theorem}

A ``resonance'' here is defined to take place iff $W(0)=0$ where
$W(\lambda)$ is the Wronskian of the two Jost solutions at energy $\lambda^2$, see the 
following section. It is known that the spectrum of $H$
is purely absolutely continuous on $(0,\infty)$ under our assumptions ($V\in L^1(\R)$ suffices for
that) so that $P_{ac}$
is the same as the projection onto the orthogonal complement of the bound states.
For the case of three dimensions we prove the following result.

\begin{theorem}
\label{thm:thm2}
Let $|V(x)|\le C(1+|x|)^{-\beta}$ for all $x\in\R^3$ where $\beta>3$.
Assume also that zero is neither an eigenvalue nor a resonance of $H=-\Laplace+V$. 
Then
\begin{equation}
\label{eq:dis_3} \big\|e^{itH} P_{ac}(H)\big\|_{1\to\infty} \les |t|^{-\frac32}. 
\end{equation}
\end{theorem}

See Section~3 for a discussion of resonances. In this case, too, it is 
well-known that the spectrum is purely absolutely continuous on $[0,\infty)$.
 
Such dispersive estimates have a long history.  For exponentially decaying
potentials Rauch~\cite{Rauch}  proved dispersive bounds in exponentially weighted
$L^2$-spaces. Jensen, Kato~\cite{JK} replaced exponential with polynomial decay
and obtained asymptotic expansions of $e^{itH}$ (in terms of powers of~$t$)
in the usual weighted $L^{2,\sigma}$ spaces. 
Journ\'e, Soffer, and Sogge~\cite{JSS} proved a version of Theorem~\ref{thm:thm2}
under the additional assumptions that $\beta>7$, $\hat{V}\in L^1$ and $V$ has some
additional small amount of regularity. Yajima~\cite{Y1} for the case $d=3$ proved
that the wave operators are bounded on $L^p(\R^3)$ for all $1\le p\le \infty$
assuming again that zero is neither an eigenvalue nor a resonance provided $\beta>5$
(and with similar conditions if $d>3$). 
As a consequence one obtains the $L^1\to L^\infty$ dispersive bounds. 
Our approach is very different from both~\cite{JSS} and~\cite{Y1}.
Journ\'e, Soffer, and Sogge use a time-dependent method and expand the evolution
repeatedly by means of Duhamel's formula. For large energies the smallness needed
to control the evolution $e^{itH}$ appearing on the right-hand side of such an expansion
is obtained from Kato's smoothing estimate. For small energies
they use the expansion of the resolvent around zero energy from~\cite{JK}.
Since their method relies on the integrability of $t^{-\frac{d}{2}}$ at infinity,
it can only be used in dimensions $d\ge3$ and it also requires
more regularity of $V$ ($\hat{V}\in L^1$ is a natural assumption for their
proof). Yajima~\cite{Y1} uses the stationary approach of Kato~\cite{Kato}
to bound the wave operators on $L^p$. While his result is more general
(it yields many more corollaries than just dispersive estimates), 
our approach to~\eqref{eq:dis_3} is direct and also requires less of~$V$.
The one-dimensional case was open until recently. Weder~\cite{Wed1}
proved a version of Theorem~\ref{thm:thm1} under the stronger assumption that
$\int_{-\infty}^\infty |V(x)|(1+|x|)^{\frac32+\eps}\, dx<\infty$. Later, Weder~\cite{Wed2},
and also Artbazar, Yajima~\cite{AY} established corresponding theorems
for the wave-operators. More precisely, they showed that the wave operators
are bounded on $L^p(\R)$ provided $1<p<\infty$ under similar assumptions on~$V$.
While our analysis is in some ways similar to Weder's~\cite{Wed1}, it turns out
that the high energy case can be treated more easily by means of a Born series
expansion, whereas small energies fall under the framework of the Jost solutions
as developed in the fundamental paper by Deift, Trubowitz~\cite{DT}. The latter
was also observed by Weder, but there is no need to impose any stronger condition
on~$V$ other than the one used in~\cite{DT}, i.e., $V\in L^1_1(\R)$. 

Dispersive estimates in two dimensions are unknown in this degree of
generality. Yajima~\cite{Y3} established the $L^p(\R^2)$ boundedness of
the wave operators under suitable assumptions on the decay of $V$ as well
as the behavior of the Hamiltonian at zero energy. Since his result requires that $1<p<\infty$,
it does not imply the $L^1(\R^2)\to L^\infty(\R^2)$ decay bounds for $e^{itH}P_{ac}$,
although it does imply the Strichartz estimates. We claim that our three-dimensional
argument can be adapted to two dimensions as well, since it does not require integrability
of $t^{-1}$ at infinity (unlike, say,~\cite{JSS}). Generally speaking, we expect the
argument to apply to any dimension (in $d=1$, however, we use a different strategy which
yields sharper results). For small energies we use expansions
of the perturbed resolvent around zero energy. These were unknown in $\R^2$ for some time, but
were recently obtained by Jensen, Nenciu~\cite{JenNen}, whereas dimensions three and higher
were treated by Jensen, Kato~\cite{JK} and Jensen~\cite{Jensen1},~\cite{Jensen2}. 
We plan to present appropriate versions of Theorem~\ref{thm:thm2} in  dimensions two, or four
and higher, elsewhere. 

An interesting issue in Theorems~\ref{thm:thm1} and~\ref{thm:thm2} is the question
of optimality. The decay rate of $(1+|x|)^{-2-\eps}$ appears to be a natural threshold 
for the dispersive estimates, and Theorem~\ref{thm:thm1} achieves this rate. But we
do not know at this point whether or not the statement of that theorem can also
hold under weaker assumptions -- the methods of proof used in this paper will
certainly no longer apply for more slowly decaying potentials in the case of Theorem~\ref{thm:thm1}. 
On the other hand, it is possible that the methods employed in the proof 
of Theorem~\ref{thm:thm2} do allow one to go below $\beta>3$. 
Let us remark that the weaker Strichartz estimates were shown to hold
under the condition $\beta>2$ in~\cite{RS} by a completely different argument. 

\section{The one-dimensional case}

Let $H=-\frac{d^2}{dx^2}+V$ in $\R^1$. Our goal is to prove Theorem~\ref{thm:thm1}. 
It is well-known that for $V\in L^1(\R)$, $H$ is essentially self-adjoint on the domain
\[ \big\{f\in L^2(\R)\:|\: f,f'\text{\ are a.c. and }-f''+Vf\in L^2(\R)\big\} \]
so that $e^{itH}$ is unitary. Hence~\eqref{eq:d1dec} is to be understood as the statement
\[ 
\|e^{itH}P_{ac}f\|_\infty \les |t|^{-\half}\|f\|_1 \text{\ \ for all\ \ }f\in L^1(\R)\cap L^2(\R),
\]
which then extends to all of $L^1(\R)$. 
We start with the high energy part of the argument.

\begin{lemma}
\label{lem:disp}
Let $\lambda_0 = \|V\|_1^2$ and suppose $\chi$ is a smooth cut-off such that
$\chi(\lambda)=0$ for
$\lambda\le \lambda_0$ and $\chi(\lambda)=1$ for $\lambda\ge 2\lambda_0$. 
Then
\[ \big\|e^{itH}\,\chi(H)\big\|_{1\to\infty} \les |t|^{-\half} \]
for all $t$.
\end{lemma}
\begin{proof}
In the limit $\eps\to 0+$ the one-dimensional resolvent $R_0(\lambda+i\eps):= \Bigl(-\frac{d^2}{dx^2}-(\lambda+i\eps)\Bigr)^{-1}$ has the kernel
\begin{equation}
\label{eq:R0}
 R_0(\lambda\pm i0)(x)=\frac{\pm i}{2\sqrt{\lambda}} e^{\pm i|x|\sqrt{\lambda}}. 
\end{equation}
Because of the decay of this kernel in $\lambda$, 
the resolvent $R_V(\lambda+i\eps)=(H-(\lambda+i\eps))^{-1}$ can be expanded into the Born series 
\begin{equation}
R_V(\lambda\pm i0) = \sum_{n=0}^\infty R_0(\lambda\pm i0)(-VR_0(\lambda\pm i0))^n. 
\label{eq:born}
\end{equation}
More precisely, since $\|VR_0(\lambda\pm i0)\|_{1\to 1} \le (2\sqrt{\lambda})^{-1}\,\|V\|_1 $, one has
\[ \big|\langle R_0(\lambda+i0)(VR_0(\lambda+i0))^n f,g \rangle \big| \le
 (2\sqrt{\lambda})^{-n-1} \|V\|_1^n \|f\|_1\, \|g\|_1, \]
so that \eqref{eq:born} converges provided $\lambda>\lambda_0=\|V\|_1^2$ in the following weak sense:
\begin{equation}
\label{eq:born2}
 \langle R_V(\lambda\pm i0)f,g \rangle = \sum_{n=0}^\infty \langle R_0(\lambda\pm i0)(-VR_0(\lambda\pm i0))^n f,g \rangle 
\end{equation}
for any pair of $L^1$ functions $f,g$. For such functions it is a standard fact that
\begin{equation}
\nn  R_V(\lambda-i0)g\in L^\infty(\R)
\end{equation}
provided $\lambda>0$ (this follows, for example, from the boundedness of the Jost solutions, see below). 
Therefore, the error term in any finite Born expansion, i.e., $R_V(\lambda+i0)(VR_0(\lambda+i0))^n$, 
tends to zero weakly as $n\to\infty$ provided $\lambda>\lambda_0$ since
\begin{align*}
&  \big|\langle R_V(\lambda+i0)(VR_0(\lambda+i0))^n f,g \rangle \big| \le
   \|(VR_0(\lambda+i0))^n f\|_1 \|R_V(\lambda-i0)g\|_\infty \\
& \le 
 (2\sqrt{\lambda})^{-n} \|V\|_1^n \|f\|_1\, \|R_V(\lambda-i0)g\|_\infty. 
\end{align*}

For technical reasons we introduce a truncated version $\chi_L$ of the cut-off $\chi$:
$ \chi_L(\lambda)= \chi(\lambda)\phi(\lambda/L) $
where $\phi$ is smooth, $\phi(\lambda)=1$ if $|\lambda|\le 1$, $\phi(\lambda)=0$ 
if $|\lambda|\ge 2$, and $L\ge1$. We need to show that
\begin{equation}
\sup_{L\ge1} \big|\langle e^{itH} \chi_L(H) f,g \rangle \big| \le C|t|^{-\half} \|f\|_1\|g\|_1
\end{equation}
for any pair of Schwartz functions $f,g$. 
The  absolutely continuous part of the spectral measure of $H$, which we denote by $E_{ac}(d\lambda)$, and the resolvent $R_V(\lambda+i0)$ are related by the well-known formula
\begin{equation}
\label{eq:stone}
\langle E_{ac}(d\lambda) f,g \rangle 
= \langle \frac{1}{2\pi i}[R_V(\lambda+i0)-R_V(\lambda-i0)]f,g \rangle \,d\lambda.
\end{equation}
Since $\chi_L(H) E(d\lambda)=\chi_L(H) E_{ac}(d\lambda)$ one concludes that
\begin{align*}
\big|\langle e^{itH} \chi_L(H) f,g \rangle \big|
&= \Big|(2\pi i)^{-1}\sum_{n=0}^\infty \int_{-\infty}^\infty e^{it\lambda^2}
\chi_L(\lambda^2)\lambda \langle R_0(\lambda^2+i0)
(VR_0(\lambda^2+i0))^n f,g\rangle d\lambda\;\Big|
\end{align*}
where we have first changed variables $\lambda \to\lambda^2$. 
Summation and integration may be exchanged because the Born series
converges absolutely in the $L^1_{loc}(d\lambda)$ norm, and the domain of integration
is extended to $\R$ from $[0,\infty)$ via the identity
$R_0(\lambda^2 - i0) = R_0((-\lambda)^2 + i0)$ (where $R_0(z^2)$ is interpreted as
an analytic function for $z\ne0$, see~\eqref{eq:R0}).
The kernel of $R_0(\lambda^2+i0)(VR_0(\lambda^2+i0))^n$ is given explicitly
by the formula
\[ R_0(\lambda^2+i0)(VR_0(\lambda^2+i0))^n(x,y) = 
\frac{1}{(2\lambda)^{n+1}}\int_{\R^n} \prod_{j=1}^nV(x_j)
e^{i\lambda(|x-x_1| + |y-x_n| +\sum_{k=2}^n|x_k-x_{k-1}|)} dx_1\ldots dx_n \]
Fubini's theorem permits integration in $d\lambda$
prior to all of the $dx_j$, leading to the desired bound
\begin{align}
\big|\langle e^{itH} \chi_L(H) f,g \rangle \big|
&\les \sum_{n=0}^\infty (2\sqrt{\lambda_0})^{-n} 
\sup_{a\in\R} \left| \int_{-\infty}^\infty e^{i(t\lambda^2+a\lambda)} 
\chi_L(\lambda^2)\,\lambda^{-n}\lambda_0^{n/2}\, d\lambda \right| \|V\|_1^n \, \|f\|_1 \|g\|_1
\label{eq:square} \\ 
&\le C(V_1)\, |t|^{-\half} \|f\|_1\|g\|_1. \label{eq:bingo}
\end{align}
We have used the dispersive bound for the one-dimensional Schr\"odinger
equation to estimate the oscillatory integral in~\eqref{eq:square}. Indeed, 
the quantity inside the absolute values is the solution of a one-dimensional
Schr\"odinger equation at time $t$ and position~$a$ with initial data 
$[\chi_L(\lambda^2)\,\lambda^{-n}\lambda_0^{(n+1)/2}]^{\vee}$. 
In order to pass to~\eqref{eq:bingo}, it therefore remains to check that 
\begin{equation}
\label{eq:pruf}
\sup_{n\ge0}\sup_{L\ge1} 
\big\| [\chi_L(\lambda^2)\,\lambda^{-n}\lambda_0^{n/2}]^{\vee} \big\|_1 < \infty. 
\end{equation}
For $n=0$ this reduces to 
\begin{equation}
\| [{\chi}_L(\lambda^2)]^{\vee} \|_1  \le \| L\widehat{\phi(\lambda^2)}(L\xi)\|_1(1 + \|(1-\chi(\lambda^2))^{\vee}\|_1) < \infty
\label{eq:n0}
\end{equation}
uniformly in $L$ since $1-\chi$ is compactly supported and smooth.
For general $n$ one has
\begin{align*}
& \| [\chi_L(\lambda^2)\,\lambda^{-n}]^{\vee}(\tau)\tau^2 \|_\infty  =
\| [(\chi_L(\lambda^2)\,\lambda^{-n})'']^{\vee}(\tau) \|_\infty 
 \le \| (\chi_L(\lambda^2)\,\lambda^{-n})'' \|_1 \le C(\lambda_0)\,\lambda_0^{-n/2}, 
\end{align*}
where the constant $C_0(\lambda)$ only depends on~$\lambda_0$, but not on $n$ or~$L$.
For $n\ge2$ one also has
\[
 \| [\chi_L(\lambda^2)\,\lambda^{-n}]^{\vee}(\tau) \|_\infty  \le \| \chi_L(\lambda^2)\,\lambda^{-n} \|_1 \le C(\lambda_0)\,\lambda_0^{-n/2}, 
\]
so that~\eqref{eq:pruf} follows for $n\ge2$. 
It remains to check that $\| [\chi_L(\lambda^2)\,\lambda^{-1}]^{\vee}(\tau) \|_\infty < \infty$
uniformly in~$L$. However, 
\[
 \|[\chi_L(\lambda^2)\,\lambda^{-1}]^{\vee}\|_{\infty} \le 
 \|[\chi_L(\lambda^2)]^{\vee}\|_1\|[\lambda^{-1}]^{\vee}\|_\infty < \infty
\]
uniformly in $L\ge1$ by~\eqref{eq:n0} and $[\lambda^{-1}]^{\vee}(\xi)=-i\,\sign(\xi)$.
\end{proof}

For the low energy part we use the Jost solutions $f_{\pm}(z,\cdot)$. They
are defined as solutions of
\[ -f_{\pm}''(z,x) +V(x)f_{\pm}(z,x)= z^2 f_{\pm}(z,x) \]
for $\Im z\ge0$ satisfying $|f_{\pm}(z,x)-e^{\pm izx}|\to 0$ as $x\to \pm \infty$.
In what follows we set $z=\lambda\in\R$. 
They are known to exist for $\lambda\ne0$ if $V\in L^1(\R)$. If $V\in L^1_1$, then
they also exist at $\lambda=0$. 
Denote their Wronskian by $W(\lambda)=W[f_+(\lambda,\cdot),f_{-}(\lambda,\cdot)]$.
It is well-known~\cite{DT} that $W(\lambda)\ne0$ if $\lambda\ne0$.
The Green's function has kernel
\begin{equation}
\label{eq:green}
 (H-(\lambda^2\pm i0))^{-1}(x,y)=R_V(\lambda^2\pm i0)(x,y) = \frac{f_+(\pm \lambda,y) f_{-}(\pm \lambda,x)}{W(\pm\lambda)} 
\end{equation}
for all $\lambda \not =0$ and $x<y$ (and the positions of $x,y$ reversed if $x>y$).
If $W(0)=0$ we say that zero energy is a resonance.  
Therefore,  in the non-resonant case, for any $x<y$, and any smooth, compactly
supported (say) cut-off $\chi$,  
\begin{align*}
& 2\pi i \int_0^\infty e^{it\lambda} \chi(\lambda) E_{a.c.}(d\lambda)(x,y) 
 =  \int_0^\infty e^{it\lambda^2} \lambda\chi(\lambda^2) \Big[\frac{f_+(\lambda,y)f_-(\lambda,x)}{W(\lambda)} - \frac{f_+(-\lambda,y)f_-(-\lambda,x)}{W(-\lambda)} \Big]\,d\lambda \nn \\
& =  \int_{-\infty}^\infty e^{it\lambda^2} \lambda\chi(\lambda^2) \frac{f_+(\lambda,y)f_-(\lambda,x)}{W(\lambda)}\, d\lambda. \label{eq:resol}
\end{align*}
In view of Lemma~\ref{lem:disp}, the non-resonant part of Theorem~\ref{thm:thm1} will
follow from the following low-energy lemma.

\begin{lemma}
\label{lem:low} 
Let $V\in L_1^1(\R)$ and $W(0)\ne0$. Then
\[ \sup_{x<y} \left| \int_{-\infty}^\infty e^{it\lambda^2} \frac{\lambda \chi(\lambda)}{W(\lambda)}
\,f_+(\lambda,y)f_-(\lambda,x)\, d\lambda \right| \les |t|^{-\half},\]
for any $t$. Here $\chi$ is any smooth, compactly supported cut-off.
\end{lemma}
We will distinguish the cases $x<0<y$,\  $0<x<y$, and $x<y<0$. Write 
\[ f_{\pm}(\lambda,x)=e^{\pm i\lambda x} m_{\pm}(\lambda,x) \]
so that $|m_{\pm}(\lambda,x)-1|\to 0$ as $x\to \pm \infty$. It is
known, see~\cite{DT}, that $m_{\pm}(z,x)-1$ as a function of~$z$ 
belongs to $H^2(\Compl^+)$, the Hardy space
on the upper half plane. Moreover, $m_{\pm}(\hat{\xi},x)-\delta_0(\xi)\in \meas$
relative to $\xi$ where $m_{\pm}(\hat{\xi},x)$ denotes the Fourier transform 
in the first variable alone and $\meas$ are the (complex) measures. 
By the $H^2$ property, $m_{\pm}(\hat{\xi},x)$ is supported in~$\xi\ge0$. 
A number of pointwise estimates can be made for $m_{\pm}(\hat{\xi},x)$.
Define $I(\xi) := \int_{|t|>|\xi|}|V(t)|dt$.  Then
\begin{equation} \label{eq:mL1}
\begin{aligned}
\sup_{x\ge 0} |m_{+}(\hat{\xi},x)-\delta_0| &\les I(\xi), &
\sup_{x\le 0} |m_{-}(\hat{\xi},x)-\delta_0| &\les I(\xi), \\ 
x\ge0\impl |\partial_x\,m_{+}(\hat{\xi},x)| &\les I(\xi) + |V(x+\xi)|,& 
|\partial_\xi(m_+(\hat{\xi},x)-\delta_0)| &\les 
 I(\xi) + |V(x+\xi)|,& \quad \\
x\le0\impl|\partial_x\,m_{-}(\hat{\xi},x)| &\les I(\xi) + |V(x-\xi)|, & 
|\partial_\xi(m_-(\hat{\xi},x)-\delta_0)| &\les
 I(\xi) + |V(x-\xi)|,
\end{aligned}
\end{equation}
see Lemma~3 in~\cite{DT}.  
The assumption $V \in L^1_1(\R)$ suffices to bound the total variation norms
$\norm[m_+(\,\hat{\cdot},x)][\meas]$, 
$\norm[\partial_xm_+(\,\hat{\cdot},x)][\meas]$ , and 
$\norm[\partial_\xi(m_+(\,\hat{\cdot},x)-\delta_0)][\meas]$
uniformly in $x \ge 0$.
Similarly, the norms $\norm[m_-(\,\hat{\cdot},x)][\meas]$, 
$\norm[\partial_xm_-(\,\hat{\cdot},x)][\meas]$, 
and $\norm[\partial_\xi(m_-(\,\hat{\cdot},x) - \delta_0)][\meas]$ 
are bounded uniformly in $x \le0$.  Identical bounds are also true of 
$[\chi(\cdot)m_{\pm}(\cdot,x)]^{\wedge}(\xi)$, however the convolution with
$\hat{\chi}$ provides enough smoothing that the norms may be taken in
$L^1(\xi)$, and the point-mass correction $\delta_0$ is not needed.

If $V$ satisfies the stronger hypothesis $V\in L^1_2(\R)$, then $I\in L^1_1(\R)$
leading to uniform bounds for $m_\pm(\hat{\xi},x)$ and its derivatives in
the $L^1_1(\R)$ norm.
Note, however, that these suprema are typically not finite if they
are taken over all $x\in\R$ rather than on the appropriate half-line.

\begin{lemma}
\label{lem:Wronskians}
Let $V\in L_j^1(\R),\ j = 1,2$, and $\tilde{\chi}$ a smooth,
compactly supported cut-off which
is identically 1 on the support of $\chi$. Then the functions
$\tilde{\chi}(\lambda)W(\lambda)$ and 
$W[f_+(\lambda,\cdot), f_-(-\lambda,\cdot)]$ both have
Fourier transform in $L_{(j-1)}^1(\R)$.
\end{lemma}
\begin{proof}
By definition,
\begin{align*}
&\tilde{\chi}(\lambda)W(\lambda) = 
 \tilde{\chi}(\lambda)\big(m_+(\lambda,0)\partial_x m_-(\lambda,0) -
   \partial_x m_+(\lambda,0)m_-(\lambda,0)\big) - 
     2i\lambda\tilde{\chi}(\lambda)m_+(\lambda,0)m_-(\lambda,0)  \\
&W[f_+(\lambda, \cdot), f_-(-\lambda, \cdot)] =
 m_+(\lambda,0)\partial_x m_-(-\lambda,0) -
  \partial_x m_+(\lambda,0)m_-(-\lambda,0)
\end{align*}  
The estimates in \eqref{eq:mL1} suffice to prove the lemma, since the Fourier
transform of each product will be a convolution of functions in 
$L_{(j-1)}^1(\R)$.
\end{proof}
\begin{proof}[Proof of Lemma \ref{lem:low}] 
In the case $x<0<y$,
\begin{align}
& \sup_{x<0<y} \left|\int_{-\infty}^\infty e^{it\lambda^2} 
\frac{\lambda \chi(\lambda)}{W(\lambda)}
f_{+}(\lambda,y) f_{-}(\lambda,x)\, d\lambda \right| \nn \\
& = \sup_{x<0<y} \left|\int_{-\infty}^\infty e^{it\lambda^2} e^{i\lambda(x-y)} \frac{\lambda \chi(\lambda)}{\tilde{\chi}(\lambda)W(\lambda)}
m_{+}(\lambda,y) m_{-}(\lambda,x)\, d\lambda \right| \les |t|^{-\half}. \label{eq:mm}
\end{align}
The final inequality follows again from the dispersive bound for the 
one-dimensional Schr\"odinger equation provided 
\begin{equation}
\label{eq:wiener}
\sup_{x<0<y} \left\| \left[\frac{\lambda \chi(\lambda)}{\tilde{\chi}(\lambda)W(\lambda)} m_{+}(\lambda,y) m_{-}(\lambda,x) \right]^\vee \right\| < \infty
\end{equation}
where the Fourier transform is with respect to $\lambda$ alone, and the
norm is in the sense of measures. Since $(\tilde{\chi} W)^\vee \in L^1(\R)$
by Lemma~\ref{lem:Wronskians}, Wiener's
lemma (see~\cite{Katz} chapter~VIII Lemma~6.3) implies that $\frac{\lambda\chi(\lambda)}{\tilde{\chi}(\lambda)W(\lambda)}$ is the Fourier
transform of an~$L^1(\R)$ function, and therefore~\eqref{eq:wiener} holds. 

By symmetry, it suffices to check the remaining case $0\le x<y$.
The danger here is that 
$f_-(0,x)$ can grow as $x\to \infty$. To deal with this issue, we expand
\[ f_-(\lambda, x) = \alpha_-(\lambda) f_+(\lambda,x)+\beta(\lambda) f_+(-\lambda,x)\]
where $\beta(\lambda)=\frac{W(\lambda)}{-2i\lambda}$ and 
$\alpha_-(\lambda)= \frac{-1}{2i\lambda} W[f_-(\lambda,\cdot),f_{+}(-\lambda,\cdot)]$. 
Hence
\begin{align*}
& \sup_{0\le x<y} \left| \int_{-\infty}^\infty e^{it\lambda^2} \frac{\lambda \chi(\lambda)}{\tilde{\chi}(\lambda)W(\lambda)} f_{+}(\lambda,y)
f_{-}(\lambda,x) \, d\lambda \right| \\
& \les  \sup_{0\le x<y} \left| \int_{-\infty}^\infty e^{it\lambda^2} e^{i\lambda(x+y)} \frac{\lambda\alpha_-(\lambda)}{\tilde{\chi}(\lambda)W(\lambda)} \chi(\lambda)
m_{+}(\lambda,y) m_{+}(\lambda,x) \, d\lambda \right| \\
& + \sup_{0\le x<y} \left| \int_{-\infty}^\infty e^{it\lambda^2} e^{i\lambda(y-x)}  \chi(\lambda)  m_{+}(\lambda,y) m_{+}(-\lambda,x) \, d\lambda \right| \les |t|^{-\half},
\end{align*}
where the final inequality again follows by noting that $\lambda \alpha_-(\lambda)$
has Fourier transform in $L^1(\R)$ and invoking the Wiener algebra.
\end{proof}

In the resonant case of Theorem~\ref{thm:thm1}, $W(0)$ vanishes and we
cannot apply the Wiener lemma to expressions with $W(\lambda)$ in the 
denominator.  
The additional hypothesis $V \in L^1_2(\R)$ ensures that 
$\frac{W(\lambda)}{\lambda}$
is continuous and nonzero everywhere (see Theorem~1 in~\cite{DT}).
The Fourier transform of $\tilde{\chi}(\lambda)W(\lambda)$ is in $L^1_1(\R)$ 
by Lemma~\ref{lem:Wronskians}, and furthermore by the identity
\[ \big(\frac{\tilde{\chi}W}{\lambda}\big)^\vee(\xi) = i \left[\int_\xi^\infty
(\tilde{\chi}W)^\vee(\eta)d\eta - \int_{-\infty}^\xi (\tilde{\chi}W)^\vee(\eta)d\eta \right]  \]
it follows that $(\tilde{\chi}\beta)^\vee \in L^1(\R)$ 
(keeping in mind that $\int_{-\infty}^\infty (\tilde{\chi}W)^\vee(\eta)d\eta=0$). 
A similar argument shows that $(\alpha_-)^\vee \in L^1(\R)$.
By rewriting every fraction with denominator $\tilde{\chi}(\lambda) W(\lambda)$
to have instead a denominator of $\tilde{\chi}(\lambda)\beta(\lambda)$
(this is done by canceling a common factor of $\lambda$ in the numerator and
denominator), the Wiener lemma may be applied precisely as above.

\section{The three-dimensional case}

Let $H=-\Lap+V$ in $\R^3$. Our goal is to prove Theorem~\ref{thm:thm2}.
We first recall the definition of a resonance. The meaning of this notion
will become clear later.

\begin{defi}
\label{def:res} 
As usual, we say that a resonance occurs at zero, provided there is a distributional
solution $f$ of the equation $(-\Laplace+V)f=0$ where for every $\sigma<-\half$ one has
$f\in L^{2,\sigma}(\R^3)\setminus L^2(\R^3)$. 
\end{defi}

Let $\chi$ be a smooth, even, cut-off function
on the line that is equal to one on a neighborhood of the
origin. Then, with $R_0(z):=(-\Lap-z)^{-1}$ and $R_V(z):=(H-z)^{-1}$, 
we need to prove that 
\begin{align}
\sup_{L\ge1}\Big|\Big \la e^{itH} \chi(\sqrt{H}/L)\Pac f,g \Big\ra\Big| &=  \sup_{L\ge1}\Big|\int_0^\infty 
e^{it\lambda^2}\lambda\, \chi(\lambda/L) \Big \la [R_V(\lambda^2+i0)-R_V(\lambda^2-i0)]f,g \Big\ra
\, \frac{d\lambda }{\pi i}\Big| \label{eq:spec_theo} \\
&\les |t|^{-\frac32}\|f\|_1\|g\|_1, \nn
\end{align}
see~\eqref{eq:stone}.
Iterating the resolvent identity yields the finite Born series
\begin{align}
R_V(\lambda^2\pm i0) &=  \sum_{\ell=0}^{2m+1} R_0(\lambda^2\pm i0)(-VR_0(\lambda^2\pm i0))^\ell \nn \\
& + R_0(\lambda^2\pm i0)(VR_0(\lambda^2\pm i0))^mVR_V(\lambda^2\pm i0)V(R_0(\lambda^2\pm i0)V)^mR_0(\lambda^2\pm i0).
\label{eq:res_ident} 
\end{align} 
Here $m$ is any positive integer. 
One needs to distinguish small $0<\lambda<\lambda_0$ from 
$\lambda>\lambda_0$, where $\lambda_0>0$ is a small
constant that will be determined by the small energy considerations below.
In the latter case, use the limiting absorption principle. In the former
case, one expands the resolvent $R_V$ around zero energy as in Jensen, Kato~\cite{JK}.
This requires assuming that zero is neither an eigenvalue nor a resonance. 
We will, however, not rely on~\cite{JK} but rederive the
expansion of the resolvent ourselves in the form needed here. 

\subsection{Large energies}

We now turn to the large energy estimates which will yield the desired bound on
\begin{align}
& \Big \la e^{itH} \chi\big(\frac{\sqrt{H}}{L}\big)\Big[1-\chi\big(\frac{\sqrt{H}}{\lambda_0}\big)\Big]\,\Pac \, f,g \Big\ra 
\nn \\
& =  \int_0^\infty 
e^{it\lambda^2}\lambda\, \chi\big(\frac{\lambda}{L}\big)\Big[1-\chi\big(\frac{\lambda}{\lambda_0}\big)\Big]\, \Big \la [R_V(\lambda^2+i0)-R_V(\lambda^2-i0)]f,g \Big\ra
\, \frac{d\lambda }{\pi i}, \label{eq:spec_theo_high} 
\end{align} 
cf.~\eqref{eq:stone}. 
Insert the resolvent expansion \eqref{eq:res_ident} into \eqref{eq:spec_theo_high}. 
The first $2m+2$ terms which do not contain the resolvent $R_V$ 
are treated as in~\cite{RS}, Section~2. This only requires that
\begin{equation}
\label{eq:katonorm}
 \|V\|_{\kato}:=\sup_{x\in\R^3} \int \frac{|V(y)|}{|x-y|}\, dy < \infty. 
\end{equation}
In particular, if $|V(x)|\les (1+|x|)^{-2-\eps}$, then this condition is satisfied.
The method from \cite{RS} gives an $L^1(\R^3)\to L^\infty(\R^3)$ bound with decay $|t|^{-\frac32}$
for those terms. For the convenience of the reader we recall the relevant parts
from~\cite{RS}. The contribution by the $(k+1)$-st term in the Born series~\eqref{eq:res_ident} is
equal to
\begin{align*}  & \int_0^\infty e^{it\lambda}\;\psi(\sqrt{\lambda}/L) 
(1-\chi(\sqrt{\lambda}/\lambda_0))\,
  \langle R_0(\lambda+i0)(VR_0(\lambda+i0))^k\,f,g \rangle \, d\lambda \\
& -\int_0^\infty e^{it\lambda}\;\psi(\sqrt{\lambda}/L) 
(1-\chi(\sqrt{\lambda}/\lambda_0))\,
  \langle R_0(\lambda-i0)(VR_0(\lambda-i0))^k\,f,g \rangle \, d\lambda
\end{align*}
which is
controlled by 
\begin{align}
& \sup_{L\ge1} \Bigl| \int_0^\infty e^{it\lambda}\;\psi(\sqrt{\lambda}/L) 
(1-\chi(\sqrt{\lambda}/\lambda_0))\,
 \Im \langle R_0(\lambda+i0)(VR_0(\lambda+i0))^k\,f,g \rangle \,
d\lambda \Bigr| \nonumber\\
& \le\int_{\R^6} |f(x_0)||g(x_{k+1})|
\int_{\R^{3k}} \frac{\prod_{j=1}^k |V(x_j)|}{\prod_{j=0}^k 4\pi |x_j-x_{j+1}|}\cdot\nonumber\\
& \qquad\qquad\qquad \cdot \sup_{L\ge1} \Bigl| \int_0^\infty e^{it\lambda}\; 
\psi(\sqrt{\lambda}/L) (1-\chi(\sqrt{\lambda}/\lambda_0))\,\sin\Bigl(\sqrt{\lambda}\sum_{\ell=0}^k |x_\ell-x_{\ell+1}|\Bigr)\,
d\lambda \Bigr|
\; d(x_1,\ldots,x_k)\,dx_0\,dx_{k+1} \label{eq:gross}\\
&\le Ct^{-\frac32}  \int_{\R^6} |f(x_0)||g(x_{k+1})| 
\int_{\R^{3k}} \frac{\prod_{j=1}^k |V(x_j)|}{(4\pi)^{k+1}\prod_{j=0}^k|x_j-x_{j+1}|}\sum_{\ell=0}^k |x_\ell-x_{\ell+1}|
\; d(x_1,\ldots,x_k)\;dx_0\,dx_{k+1} \label{eq:lem1} \\
&\le Ct^{-\frac32}  \int_{\R^6} |f(x_0)||g(x_{k+1})|\; 
(k+1) (\|V\|_{\kato}/4\pi)^k \;dx_0\,dx_{k+1} \label{eq:lem2} \\
&\le C_k\,t^{-\frac32} \|f\|_1\|g\|_1. \nonumber
\end{align}
In order to pass to \eqref{eq:gross} one uses the explicit representation 
of the kernel of $R_0(\lambda+i0)(x,y)=\frac{e^{ i\lambda|\x-\y|}}{4\pi|\x-\y|}$, 
which leads to a $k$-fold integral. The inequalities~\eqref{eq:lem1} and~\eqref{eq:lem2}
are obtained by means of the following two lemmas from~\cite{RS}. We provide the proof
of the first lemma, as its statement differs slightly from the one in~\cite{RS} (by the
introduction of an additional zero energy cut-off).

\begin{lemma}
\label{lem:statphas}
Let $\psi$ be a smooth, even  bump function
with $\psi(\lambda)=1$ for $-1\le\lambda\le 1$ and $\supp(\psi)\subset[-2,2]$. 
Then for all $t\ge1$ and any real~$a$,
\begin{equation}
\label{eq:decay}
\sup_{L\ge 1}\Bigl| \int_0^\infty e^{it\lambda} \sin(a\sqrt{\lambda})\,
\psi\Bigl(\frac{\sqrt\lambda}{L}\Bigr)\,(1-\chi(\sqrt{\lambda}/\lambda_0))
\,d\lambda\Bigr| \le
C \,t^{-\frac32}\,|a|
\end{equation}
where $C$ only depends on ~$\psi$, $\chi$, and $\lambda_0$.
\end{lemma}
\begin{proof}
Denote the integral in \eqref{eq:decay} by $I_L(a,t)$. 
The change of variables $\lambda\to \lambda^2$ leads to the expression
$$
I_L(a,t) = 2 \int_0^\infty \lambda \,e^{it\lambda^2} \sin (a\lambda)\,\psi(\lambda/L)\,(1-\chi(\lambda/\lambda_0))
\,d\lambda 
$$
Integrating by parts we obtain
$$
I_L(a,t) = -\frac i{t}\int_0^\infty e^{it\lambda^2} 
\bigg(a\,\cos (a\lambda)\,\psi(\lambda/L)\,(1-\chi(\lambda/\lambda_0))
 + \sin (a\lambda)\,\big[(1-\chi(\lambda/\lambda_0))\,\psi(\lambda/L)\big]'\bigg) 
\,d\lambda. 
$$
Since $\psi$ and $\chi$ are  assumed to be even, the derivative of the brackets is odd. Hence, 
$$
\aligned
I_L(a,t)=&  -\frac i{2t}\int_{-\infty}^\infty e^{it\lambda^2} 
\bigg(a\,\cos (a\lambda)\,(1-\chi(\lambda/\lambda_0))
\psi(\lambda/L) + \sin (a\lambda)\,\big[(1-\chi(\lambda/\lambda_0))\,\psi(\lambda/L)\big]'\bigg) 
\,d\lambda \\ =&  
-\frac {a}{4t}\,i\int_{-\infty}^\infty e^{it\lambda^2} 
\big(\,e^{i a\lambda} + e^{-i a\lambda}\big)\,\psi(\lambda/L)\,(1-\chi(\lambda/\lambda_0))
\,d\lambda  \\
& + \int_0^a  \frac {1}{4t}\,\int_{-\infty}^\infty e^{it\lambda^2} 
\big(\,e^{i b\lambda} + e^{-i b\lambda}\big)\,{\lambda}\big[(1-\chi(\lambda/\lambda_0))\,\psi(\lambda/L)\big]'\,d\lambda \,db.
\endaligned
$$ 
Invoking the $|t|^{-\half}$ dispersive bound for the one-dimensional Schr\"odinger
equation, it thus  suffices to show that 
\[ \|[(1-\chi(\lambda/\lambda_0))\,
\psi(\lambda/L)]^{\vee}\|_{1}+\Big\|\Big[\lambda\big[(1-\chi(\lambda/\lambda_0))\,\psi(\lambda/L)\big]'\Big]^{\vee}\Big\|_{1} < \infty
\]
uniformly in $L\ge1$. These properties are elementary and left to the reader.
\end{proof}

The following lemma is identical with one in Section 2 of~\cite{RS}, and we refer
the reader to that paper for the simple proof.

\begin{lemma}
\label{lem:iter} 
For any positive integer $k$ and $V$ as in \eqref{eq:katonorm}
\begin{equation}
\nn
\sup_{x_0,x_{k+1}\in\R^3}\int_{\R^{3k}} \frac{\prod_{j=1}^k |V(x_j)|}{\prod_{j=0}^k|x_j-x_{j+1}|}\sum_{\ell=0}^k |x_\ell-x_{\ell+1}|\; dx_1\ldots\,dx_k \le (k+1) \|V\|_{\kato}^k.
\end{equation}
\end{lemma}

We now turn to the term in the Born series~\eqref{eq:res_ident} containing the
perturbed resolvent~$R_V$. 
Recall from Agmon~\cite{agmon} or Reed, Simon~\cite{RS4} that (for general dimensions $\R^d$)
\[ \|(-\Lap -(\lambda^2\pm i0))^{-1} f\|_{L^{2,-\sigma}} \les \|f\|_{L^{2,\sigma}},\]
provided $\sigma>\half$ and $1<\lambda<2$, say. This bound is known as the limiting
absorption principle. 
It  extends easily to $\lambda>1$ with a constant that decays like
 $\lambda^{-2+2\sigma}$ (this is not optimal but sufficient for our purposes). 
Indeed, use that $\lambda^{d-2}(-\Lap-1)^{-1}(\lambda x)
=(-\Lap-\lambda^2)^{-1}(x)$ for the kernels of the resolvents.   
Since $L^{2,\alpha}$ embeds in $L^{2,\sigma}$ for all $\alpha > \sigma$, it is
to our advantage to choose $\sigma = \frac12 + $ so that
\begin{equation}
\label{eq:r0dec}
\|(-\triangle-(\lambda^2 \pm i0))^{-1}f\|_{L^{2,-\sigma}} \les
 \lambda^{-1+} \|f\|_{L^{2,\sigma}}
\end{equation}
for all $\sigma > \frac12$, $\lambda > 1$. Throughout this paper, the notation
$a+$ or $a^{+}$ for some number $a$ means $a+\eps$ for an arbitrarily small, but fixed $\eps>0$.
Similarly with $a-$ and $a^{-}$. 
The free resolvent $R_0(\lambda^2\pm i0):= (-\Lap-(\lambda^2\pm i0))^{-1}$
satisfies the following well-known bounds.  

\begin{proposition}
\label{prop:r0dev}
The derivatives $\frac{d^j}{d\lambda^j}\big[R_0(\lambda^2\pm i0)\big]$
satisfy the uniform bounds
$$\sup_{\lambda}\left\|\frac{d^j}{d\lambda^j}\big[R_0(\lambda^2\pm i0)\big]f \right\|_{
L^{2,-\sigma}} \les \|f\|_{L^{2,\sigma}}$$
for all $\sigma > j + \frac12$ and $j\ge1$. 
\end{proposition}
\begin{proof}
The kernel of $\frac{d^j}{d\lambda^j}\big[R_0(\lambda^2\pm i0)\big]$ has
the explicit form
$$\frac{d^j}{d\lambda^j}\big[R_0(\lambda^2\pm i0)\big](x,y) = \frac1{4\pi}
e^{\pm i\lambda|x-y|}|x-y|^{j-1}$$
The Hilbert-Schmidt norm of this operator as a mapping from $L^{2,\sigma}$ to
$L^{2,-\sigma}$ is given by
$$\|\frac{d^j}{d\lambda^j}\big[R_0(\lambda^2\pm i0)\big]\|_{HS}^2
 = C\iint_{\R^6}\langle x\rangle^{-2\sigma}|x-y|^{2j-2}
  \langle y\rangle^{-2\sigma} dx\,dy$$
The integral may be divided into the three domains $|x| \le \frac{|y|}2$,
$|x-y| \le \frac{|y|}2$, and the complement of these two.  For a fixed point 
$y \in \R^3$, the respective regions contribute 
$\langle y\rangle^{2j-2-2\sigma}$, $\langle y\rangle^{2j+1-4\sigma}$, and
$\langle y\rangle^{2j+1-4\sigma}$ again when integrated with respect to $x$.
If $\sigma > j + \frac12$, each of these exponents is less than $-3$, leading
to a convergent integral in $dy$.  Note that all dependence on $\lambda$
was removed by taking absolute values.
\end{proof}

Next, one transfers these estimates to $R_V(\lambda^2\pm i0)$ by
means of the resolvent identity
\begin{align}
 R_V(\lambda^2\pm i0) &= R_0(\lambda^2\pm i0) - R_0(\lambda^2\pm i0)VR_V(\lambda^2\pm i0)  \nn \\
 R_V(\lambda^2\pm i0) &= (I+R_0(\lambda^2\pm i0)V)^{-1} R_0(\lambda^2\pm i0).  \nn 
\end{align}
Now $S=S(\lambda) := I+R_0(\lambda^2\pm i0)V$ is a perturbation of the identity by the
compact operator $R_0(\lambda^2\pm i0)V: L^{2,-\sigma}\to L^{2,-\sigma}$ with $\sigma>\half$
provided $|V(x)|\les (1+|x|)^{-1-}$. The compactness here follows
from the fact that the resolvent gains two derivatives in the weighted $L^2$ space.
Thus $S^{-1}$ exists iff $Sf=0$ implies $f=0$ for any $f\in L^{2,-\sigma}$. But $Sf=0$ is
formally equivalent to $(-\Lap+V)f=\lambda^2 f$. Since $\lambda>0$, it follows from Agmon~\cite{agmon}
that in fact $f$ which was only assumed to be in $L^{2,\sigma}$ for every $\sigma>\half$, 
has to be an eigenfunction (i.e., in $L^2$). But positive embedded eigenvalues do not exist 
by Kato's theorem, see~\cite{RS4}, Section~XIII.8 for all this.  Hence $S(\lambda)^{-1}:L^{2,-\sigma}
\to L^{2,-\sigma}$ exists for all $\lambda>0$ provided $\sigma>\half$.
Furthermore, $S(\lambda)$ converges to the identity operator as 
$\lambda\to\infty$ which then implies that $S(\lambda)^{-1}$ is uniformly bounded
for all $\lambda > \lambda_0$.  Consequently, for $\sigma=\half+$,
\begin{equation}
\label{eq:rvdec}
 \|R_V(\lambda^2\pm i0)\|_{L^{2,\sigma}\to L^{2,-\sigma}} \les \lambda^{-1+}.
\end{equation}
To handle derivatives of $R_V(\lambda^2\pm i0)$, one checks that
\begin{equation}
\label{eq:der_id}
 \frac{d}{d\lambda}R_V(\lambda^2\pm i0) = -S(\lambda)^{-1}\frac{d}{d\lambda}R_0(\lambda^2\pm i0)\,VS(\lambda)^{-1}
R_0(\lambda^2\pm i0)+ S(\lambda)^{-1}\frac{d}{d\lambda} R_0(\lambda^2\pm i0),
\end{equation}
and since $\sup_{\lambda>\lambda_0}\|S(\lambda)^{-1}\|_{L^{2,-\sigma}\to L^{2,-\sigma}} < \infty$
for $\sigma>\half$, it follows that also
\begin{equation}
\label{eq:rvdev} \sup_{\lambda>\lambda_0}\|\frac{d}{d\lambda}\,R_V(\lambda^2\pm i0)\|_{L^{2,\sigma}\to L^{2,-\sigma}} \les 1 
\text{\ \ for\ \ } \sigma>\frac32.
\end{equation} 
Note from~\eqref{eq:der_id} that one needs to assume the decay $|V(x)|\les (1+|x|)^{-2-\eps}$ for
this to hold. Indeed, $V$ needs to take $L^{2,-\half-} \to L^{2,\frac32+}$.
By a similar argument,
\[ \|\frac{d^2}{d\lambda^2} R_V(\lambda^2\pm i0)\|_{L^{2,\sigma}\to L^{2,-\sigma}} \les 1
\text{\ \ for\ \ }\sigma>\frac52. \]
This estimate requires the decay $|V(x)|\les (1+|x|)^{-3-}$ by an analogous formula
to~\eqref{eq:der_id}. 

Let $R_0^{\pm}(\lambda^2):= R_0(\lambda^2\pm i0)$.
Moreover, set
\[ G_{\pm,x}(\lambda^2)(x_1):= e^{\mp i\lambda|x|}R_0(\lambda^2\pm i0)(x_1,x) 
= \frac{e^{\pm i\lambda(|x_1-x|-|x|)}}{4\pi|x_1-x|}. \]
Similar kernels appear already in Yajima's work~\cite{Y3} (see his high
energy section). 
Removing $f,g$ from~\eqref{eq:spec_theo}, we are led to proving that 
\begin{align}
& \left| 
\int_0^\infty e^{it\lambda^2}e^{\pm i\lambda(|x|+|y|)}\;\chi(\lambda/L)\, (1-\chi(\lambda/\lambda_0))
\lambda \Big\la VR^{\pm}_V(\lambda^2)V (R_0^{\pm}(\lambda^2)V)^m G_{\pm,y}(\lambda^2), (R_0^{\mp}(\lambda^2)V)^m G_{\pm,x}^*(\lambda^2) \Big\ra \, d\lambda 
\right|  \label{eq:main} \\
& \les |t|^{-\frac32} \nn
\end{align}
uniformly in $x,y\in\R^3$ and $L\ge 1$. 
\begin{proposition}
\label{prop:Gest}
The derivatives of $G_{+,x}(\lambda^2)$ satisfy the 
estimates
\begin{equation}
\label{eq:Gest}
\begin{aligned}
\sup_{x\in\R^3} \Big\| \frac{d^j}{d\lambda^j} G_{+,x}(\lambda^2)\Big\|_{L^{2,-\sigma}} &< C_{j,\sigma} \text{\ \ provided\ \ } \sigma > \frac12 + j \\
\sup_{x\in\R^3} \Big\| \frac{d^j}{d\lambda^j} G_{+,x}(\lambda^2)\Big\|_{L^{2,-\sigma}} &< \frac{C_{j,\sigma}}{\la x\ra}
\text{\ \ provided\ \ } \sigma>\frac32+j
\end{aligned}
\end{equation}
for all $j\ge 0$. 
\end{proposition}
\begin{proof}
This follows from the explicit formula
\begin{align*}
\left\|\frac{d^j}{d\lambda^j} \frac{e^{i\lambda(|u-x|-|x|)}}{|x-u|} \la u\ra^{-\sigma} \, du \right\|_2
&=  \left( \int_{\R^3} \frac{(|u-x|-|x|)^{2j}}{|x-u|^2}\, \la u\ra^{-2\sigma}\, du\right)^{\half} \\
&\le  \left( \int_{\R^3} \frac{\la u\ra^{2(j-\sigma)}}{|x-u|^2}\, du\right)^{\half} 
\end{align*}

The final estimate on this integral is obtained by dividing $\R^3$ into the regions $|u|<\frac{|x|}2$, 
$|x-u|<\frac{|x|}2$, and the complement of these two.  If $\frac12 < (\sigma - j) <
\frac32$, then each of these regions contributes $\la x\ra^{\frac12+j-\sigma}$
to the total.
If $\sigma > \frac32 + j$, the first region instead contributes
$\la x \ra^{-1}$, making it the dominant term.
\end{proof}
 
Rewrite the integral in~\eqref{eq:main} in the form (with $L=\infty$)
\begin{equation}
\label{eq:Ipm}
I^{\pm}(t,x,y):=\int_0^\infty e^{it\lambda^2 \pm i\lambda(|x|+|y|)} a^{\pm}_{x,y}(\lambda)\, d\lambda.
\end{equation}
Then in view of \eqref{eq:r0dec}, \eqref{eq:rvdec}, \eqref{eq:rvdev}, and Propositions~\ref{prop:r0dev} and~\ref{prop:Gest}, 
one concludes that $a^{\pm}_{x,y}(\lambda)$ has two
derivatives in~$\lambda$ and  
\begin{equation}
\label{eq:adec}
\begin{aligned}
\Big|\frac{d^j}{d\lambda^j} a^{\pm}_{x,y}(\lambda)\Big| &\les (1+\lambda)^{-2^+} (\la x\ra\la y\ra)^{-1} 
\text{\ \ for\ \ } j = 0,1, \text{\ \ and all\ \ }\lambda>1  \\
\Big|\frac{d^2}{d\lambda^2} a^{\pm}_{x,y}(\lambda)\Big| &\les (1+\lambda)^{-2^+} \qquad \text{\ \ for\ all\ \ }\lambda>1,
\end{aligned}
\end{equation}
which in particular justifies taking $L=\infty$ in~\eqref{eq:Ipm}.
This requires that one takes $m$ sufficiently large ($m=2$ is sufficient) and
that $|V(x)|\les (1+|x|)^{-\beta}$ for some $\beta>3$. The latter condition 
arises as follows: Consider, for example, the case where two derivatives fall
one of the $G$-terms at the ends. Then $V$ has to compensate for $\frac52^+$ powers
because of~\eqref{eq:Gest}, and also a $\frac12^+$ power from
\[ \|R_0^{\pm}(\lambda^2) f\|_{L^{2,-\half-}} \les \lambda^{-1^+} \|f\|_{L^{2,\half+}}. \]
Similarly with the other terms. 

As far as $I^{+}(t,x,y)$ is concerned, note that 
on the support of~$a^{\pm}_{x,y}(\lambda)$ the phase $t\lambda^2+\lambda(|x|+|y|)$ has 
no critical point.  Two integrations by parts yield the bound 
$|I^{+}(t,x,y)|\les t^{-2}$. 

In the case of $I^{-}(t,x,y)$ the phase $t\lambda^2-\lambda(|x|+|y|)$
has a unique critical point at $\lambda_1=(|x|+|y|)/(2t)$. If $\lambda_1\ll \lambda_0$, then
two integration by parts again yield a bound of~$t^{-2}$.
If $\lambda_1\gtrsim \lambda_0$ then the bound
$\max(|x|,|y|)\gtrsim t$ is also true, and stationary phase contributes
$t^{-\half}(\la x\ra \la y\ra)^{-1}\les t^{-\frac32}$,
as desired.  Strictly speaking, these estimates are only useful when $t>1$. 
On the other hand, when $0<t<1$ there is nothing to prove since
$I^{\pm}(t,x,y) \les 1$ by~\eqref{eq:adec}.

To apply stationary phase properly, one should
restrict $a^{\pm}_{x,y}(\lambda)$ to a compact interval of the form $[\lambda_1-C,\lambda_1+C]$
for some constant $C\gg1$. Outside of this interval, one uses the decay given
by~\eqref{eq:adec} in terms of~$\lambda$. Two integrations by parts yield the bound
$t^{-3}$ for the remaining piece of~$I^{-}(t,x,y)$. This concludes the 
high-energy part of the argument.

\subsection{Low energies}

In view of~\eqref{eq:spec_theo} and~\eqref{eq:spec_theo_high} it remains
to control the low-energy part
\begin{align}
& \Big \la e^{itH} \chi(\sqrt{H}/\lambda_0)\,\Pac \, f,g \Big\ra 
\nn \\
& =  \int_0^\infty 
e^{it\lambda^2}\lambda\, \chi(\lambda/\lambda_0)\, \Big \la [R_V(\lambda^2+i0)-R_V(\lambda^2-i0)]f,g \Big\ra
\, \frac{d\lambda }{\pi i} \label{eq:spec_theo_low} 
\end{align} 

If $f,g \in L^1$, this can be done by evaluating the supremum
\begin{equation}
\label{eq:sup1}
\sup_{\x, \y\in\R^3} \Big| \int_0^\infty e^{it\lambda^2} \lambda
  \chi(\lambda/\lambda_0) [R_V^+(\lambda^2) - R_V^-(\lambda^2)](\x,\y)
  d\lambda \Big|
\end{equation}

We will use the resolvent identity 
\begin{equation}
\label{eq:res_id}R_V^\pm(\lambda^2) = R_0^\pm(\lambda^2) - R_0^\pm(\lambda^2)V
 (I + R_0^\pm(\lambda^2)V)^{-1} R_0^\pm(\lambda^2)
\end{equation}

The resolvents $R_0^\pm(\lambda^2)$ have an explicit kernel representation
$$R_0^\pm(\lambda^2)(\x,\y) = \frac{e^{\pm i\lambda|\x-\y|}}{4\pi|\x-\y|}$$
The numerator of this expression always has complex magnitude 1, therefore 
the size of $|R_0^\pm(\lambda^2)|$ does not depend on $\lambda$.  We will
now estimate the Hilbert-Schmidt norm of $R_0^\pm(\lambda^2)$ as a linear
map between the weighted spaces $L^{2,\sigma}$ and $L^{2,-\alpha}$.
Let
$$\norm[R][HS(\sigma,-\alpha)]^2 = \iint_{\R^6} \langle \x\rangle^{-2\sigma}
|R(\x,\y)|^2 \langle\y\rangle^{-2\alpha}\,d\x d\y$$
denote this norm. The following proposition is a well-known bound on the 
free resolvents.

\begin{proposition} \label{prop:Rnorms}
If $\sigma, \alpha > \frac12$, and $\sigma + \alpha > 2$, then 
$$\sup_{\lambda}\norm[R_0^\pm(\lambda^2)][HS(\sigma,-\alpha)] \le C_{\sigma,\alpha}$$
\end{proposition}

\begin{proof}
The integral
$$\iint_{\R^6}\langle\x\rangle^{-2\sigma}\frac1{|\x-\y|^2}\langle\y\rangle
^{-2\alpha}d\x d\y$$
may be broken up into three disjoint domains:

Domain 1: $|\x| \le \frac12|\y|$, which requires $|\x-\y| \sim |\y|$.  The
integral over Domain 1 contributes less than 
$\int_{\R^3} \langle\y\rangle^{3-2\sigma} \langle\y\rangle^{-1-2\alpha} d\y$,
which is bounded by a constant $C_{\sigma,\alpha}$,
to the total integral.

Domain 2:  $|\x-\y| \le \frac12|\y|$, which requires $|\x|\sim|\y|$.  
The integral over Domain 2 contributes less than $\int_{\R^3} |y|
\langle\y\rangle^{-2\sigma-2\alpha} d\y$, which is also bounded by 
$C_{\sigma,\alpha}$, to the total integral.

Domain 3:  $|\x|,|\x-\y| \ge \frac12|\y|$, which requires $|\x|\sim|\x-\y|$.
The integral over Domain 3 contributes less than 
$\int_{\R^3}\langle\y\rangle^{1-2\sigma} \langle\y\rangle^{-2\alpha} d\y 
\les C_{\sigma,\alpha}$ to the total integral.
\end{proof}

If $|V(\x)| \les \langle\x\rangle^{-\beta}$ for some $\beta > 3$,
it follows that the operator
$R_0^\pm(\lambda^2)V$ is compact on the weighted space $L^{2,\sigma}(\R^3)$
for all choices of $-\frac52 \le \sigma < -\frac12$.  Indeed, one checks
by means of Proposition~\ref{prop:Rnorms} that
$R_0^\pm(\lambda^2)V$ maps $L^{2,\sigma}(\R^3)$ compactly into $L^{2,\sigma+1}
(\R^3)$ for all $\sigma \in [-\frac52,-\frac32)$.

Let $S_0 = I + R_0(0)V$. By compactness of $R_0(0)V$, the invertibility of $S_0$ depends only on whether
a solution exists in $L^{2,\sigma}$ to the equation $\psi = -R_0(0)V\psi$.
However if a solution $\psi$ satisfies  $\psi\in L^{2,\sigma}$ for some 
$\sigma \ge -\frac52$, then 
$\psi = -R_0(0)V\psi \in L^{2,\alpha}$ for any choice of $\alpha < -\frac32$.
Applying the bootstrapping process again, we see that the solution $\psi$
must lie in $L^{2,\alpha}$ for all $\alpha < -\frac12$.

It is easy to see that this same function $\psi$ is also a distributional
solution to $(\triangle +V)\psi = 0$. Conversely, any distributional solution
of $(\triangle +V)\psi = 0$ with $\psi\in L^{2,-\half-}$ satisfies $S_0 \psi=0$.  
It follows that $S_0$ is invertible in
 $L^{2,\sigma},\, -\frac52 \le \sigma < -\frac12$ precisely when zero energy
is neither an eigenvalue nor a resonance of the potential $V$, see Definition~\ref{def:res}.

Write $R_0^\pm(\lambda^2) = R_0(0) + B^\pm(\lambda)$.  Then
$$[I+R_0^\pm(\lambda^2)V]^{-1} = S_0^{-1}[I+B^\pm(\lambda)VS_0^{-1}]^{-1}$$
Examining the kernel,
$$B^\pm(\lambda)(\x,\y) = \frac{e^{\pm i\lambda|\x-\y|} - 1}{4\pi |\x-\y|}$$
which satisfies the size estimates
\begin{equation} \label{eq:Bkernel}
\big|B^\pm(\lambda)(\x,\y)\big| \ \lesssim \ \left\{
\begin{aligned}
  &\lambda &{\rm if}\ |\x-\y| \le 1/\lambda \\
  |\x&-\y|^{-1} &{\rm if}\ |\x-\y| \ge 1/\lambda
\end{aligned}
\right.
\end{equation}

The first $\lambda$-derivative of $B^\pm$ has kernel
$(B^\pm)'(\lambda)(\x,\y) =
  \frac{\pm i}{4\pi}e^{\pm i\lambda |\x-\y|}$
with the obvious bound $|(B^\pm)'(\lambda)(\x,\y)| \le C$.

The symmetry between $B^+$ and $B^-$ is expressed by the relationship
$$B^-(\lambda) = B^+(-\lambda) \quad {\rm for\ all\ }\lambda \ge 0$$

\begin{proposition} \label{prop:Bnorms1}
If $\sigma, \alpha > \frac12$, and $\sigma + \alpha > 2$, then
\ $\lim_{\lambda\to 0} \norm[B^\pm(\lambda)][HS(\sigma,-\alpha)] = 0$.
\end{proposition}

\begin{proof}
The kernels $B^\pm(\lambda)(\x,\y)$ are pointwise dominated by 
$\frac1{|\x-\y|}$, which has a finite $HS(\sigma,-\alpha)$ norm 
by Proposition~\ref{prop:Rnorms}. The result then follows from
the dominated convergence theorem.
\end{proof}

\begin{cor}
If $|V(x)| \les \langle x\rangle^{-\beta}$ for some choice of $\beta > 3$, then
$$\lim_{\lambda\to 0} \norm[B^\pm(\lambda)VS_0^{-1}][HS(\sigma,\sigma)] = 0 $$
for all $\sigma \in (-\frac52, -\frac12)$.
\end{cor}
\begin{proof} One has $VS_0^{-1}:L^{2,\sigma}\to L^{2,\sigma+3+}$
provided that $-\frac52<\sigma<-\half$. The proposition implies that
$\|B^{\pm}(\lambda)\|_{HS(\sigma+3+,\sigma)}\to0$ as $\lambda\to0$.
\end{proof}

\smallskip \noindent
{\bf Claim:}  $\norm[(\Bp)'(\lambda)][HS(\sigma,-\alpha)] \le C$\ if\ 
$\sigma,\alpha > \frac32$.
\begin{proof}
This is trivial because the function $\langle\x\rangle^{-2\sigma}
\langle\y\rangle^{-2\alpha}$
is integrable over $\R^6$.
\end{proof}

For sufficiently small $\lambda < \lambda_0$, it is then possible to expand
$$\tilde{B}^\pm(\lambda) = [I + B^\pm(\lambda)VS_0^{-1}]^{-1}$$
as a Neumann series in the norm $\norm[\cdot][HS(\sigma,\sigma)]$ for
all values $-\frac52 < \sigma < -\frac12$.

The symmetry $\tilde{B}^-(\lambda) = \Btp(-\lambda)$ is still valid.

\smallskip
For ease of notation, define $\chi_0(\lambda) = \chi(\lambda/\lambda_0)$
and $\chi_1(\lambda) = \chi(\lambda/2\lambda_0)$.  Note that $\chi_1\chi_0
= \chi_0$.
In view of \eqref{eq:sup1} and~\eqref{eq:res_id} we wish to control the size of

\begin{equation*}
\begin{aligned}
\sup_{\x,\y\in\R^3}\bigg|\int_0^\infty e^{it\lambda^2} \lambda \chi_0(\lambda)
\Big[ &\big[\Rnp - R_0^-(\lambda^2)\big] \\
 -\ &\big[\Rnp VS_0^{-1}\Btp(\lambda)\Rnp - R_0^-(\lambda^2)VS_0^{-1}
   \tilde{B}^-(\lambda)R_0^-(\lambda^2)\big]\Big](\x,\y)\, d\lambda \bigg|
\end{aligned}
\end{equation*}
\begin{equation*}
\begin{aligned}
\le\ \sup_{\x,\y\in\R^3}
\Big|\int_{-\infty}^\infty &e^{it\lambda^2}\lambda \chi_0(\lambda)
\frac{e^{i\lambda|\x-\y|}}{4\pi |\x-\y|} d\lambda \Big|
\\
+\ \sup_{\x,\y\in\R^3} \Big|&\int_{-\infty}^\infty e^{it\lambda^2}\lambda
   \iint_{\R^6} \frac{V(x_4)e^{i\lambda|\y-x_4|}}{|\y-x_4|}
  \big(S_0^{-1} (\chi_0\Btp)(\lambda)(x_4,x_1)\big)
  \frac{e^{i\lambda|\x-x_1|}}{|\x-x_1|}\, dx_1 dx_4 d\lambda \Big|
\end{aligned}
\end{equation*}

\smallskip
The first term is simply the low-energy part of the free Schr\"odinger
evolution, which is known to be dispersive.

\smallskip
The second term can be integrated by parts once, leaving
\begin{equation} \label{eq:IBP}
\sup_{\x,\y \in\R^3} \frac1{2t} \Big|
 \int_{-\infty}^\infty e^{it\lambda^2} \iint_{\R^6}\dfrac{d}{d\lambda}\Big[
 \frac{V(x_4)e^{i\lambda|\y-x_4|}}{|\y-x_4|} 
  \big(S_0^{-1} (\chi_0 \Btp)(\lambda)(x_4,x_1)\big)
 \frac{e^{i\lambda|\x-x_1|}}{|\x-x_1|} \Big] \,dx_1 dx_4 d\lambda  \Big|
\end{equation}
to be controlled.  Consider the term where $\frac{d}{d\lambda}$ falls on
$\Btp(\lambda)$.  The others will be similar.  

Using Parseval's identity, and the fact that
$\norm[(e^{it(\cdot)^2})^\wedge(u)][L^\infty(u)] = C t^{-1/2}$,
this is less than
$$\sup_{\x,\y\in\R^3} \frac1{t^{3/2}} \int_{-\infty}^\infty
  \Big|  \iint_{\R^6} \frac{V(x_4)}{|\y-x_4|} S_0^{-1}
  \big[\chi_0(\Btp)'\big]^\vee 
  \big(u+|\y-x_4| + |\x-x_1|\big)(x_4,x_1) \frac{1}{|\x-x_1|}\, dx_1 dx_4
  \Big| \, du$$

If the absolute value is taken inside the inner integral, then Fubini's
theorem may be used to exchange the order of integration to obtain

\begin{equation*}
\begin{aligned}
\sup_{\x,\y\in\R^3} \frac1{t^{3/2}} \iint_{\R^6} \int_{-\infty}^\infty
\frac{|V(x_4)|}{|\y-x_4|}  \Big| S_0^{-1}\big[\chi_0(\Btp)'\big]^\vee 
\big(u + |\y-x_4| + |\x- x_1|\big) 
{\scriptstyle (x_4, x_1)} \Big|
\frac{1}{|\x-x_1|} \, du\, dx_1 dx_4
\\
\le \sup_{\x,\y\in\R^3} \frac1{t^{3/2}} 
\Big\|\frac{|V(\cdot)|}{|\y- \cdot|}\Big\|_{L^{2,2^+}}
 \ \big\|{\textstyle \int} |S_0^{-1} [\chi_0(\Btp)']
  ^\vee(u)|du \big\|_{HS(-1^-,-2^-)} \ \big\||\x-\cdot|^{-1}\big\|_{L^{2,-1^-}}
\end{aligned}
\end{equation*}
The weighted $L^{2,-1^-}(dx_1)$-norm of $|\x - x_1|^{-1}$ is uniformly bounded
for all choices of $\x \in \R^3$.  In fact, these functions are even bounded
in the weaker $L^{2,\sigma}$norm for any $\sigma < -\frac12$.  Similarly,
the functions $\frac{V(x_4)}{|\y-x_4|}$ are uniformly bounded in
$L^{2,\sigma}(dx_4)$ for any $\sigma < \beta - \frac12$.  We are assuming
$\beta > 3$, which is more than sufficient.  It therefore remains
only to control the size of
$$ \big\|{\textstyle \int} |S_0^{-1} [\chi_0(\Btp)']^\vee(u)|du 
  \big\|_{HS(-1^-,-2^-)}$$
Minkowski's Inequality allows us to bring the norm inside the integral.
Recall that $S_0^{-1}$ is a bounded operator on $L^{2,-2^-}$, and that
the composition of a bounded operator and a Hilbert-Schmidt operator is
also Hilbert-Schmidt.
The problem then reduces to establishing existence of a number $\lambda_0 > 0$
such that
\begin{equation} \label{eq:desired}
 \int_{-\infty}^\infty
 \norm[[\chi_0(\Btp)']^\vee(u)][HS(-1^-,-2^-)]\, du  < \infty
\end{equation}

The operators $\Btp(\lambda)$ were originally defined by the convergent Neumann
series
$$\Btp(\lambda) \ = \ [I + \Bp(\lambda)VS_0^{-1}]^{-1}\  =\  \sum_{n=0}^\infty
 \big(- \Bp(\lambda)VS_0^{-1}\big)^n$$
Thus
\begin{equation} \label{eq:Bprimeseries}
\chi_0(\lambda)(\Btp)'(\lambda) \ = \ \sum_{n=1}^\infty\sum_{m=0}^{n-1} (-1)^n
\big((\chi_1\Bp)(\lambda)VS_0^{-1}\big)^m 
\chi_0(\lambda)(\Bp)'(\lambda)VS_0^{-1}
\big((\chi_1\Bp)(\lambda)VS_0^{-1}\big)^{n-(m+1)}
\end{equation}
We will take the Fourier transform of $\chi_0(\Btp)'$ term-wise and determine 
that the resulting series is convergent in the norm $L^1(du; HS(-1^-,-2^-))$.
The following refinement of Proposition~\ref{prop:Bnorms1} is especially
useful.

\begin{proposition} \label{prop:Bnorms}
Suppose $\sigma, \alpha > \frac12$, and $\alpha+\sigma > 2$.  Let $K(\lambda)$
be an integral operator on $\R^3$ whose kernel $K(\lambda)(\x,\y)$ satisfies
the size estimates in \eqref{eq:Bkernel}.  Then
$$\norm[K(\lambda)][HS(\sigma,-\alpha)] \le C_{\sigma,\alpha,\gamma}
|\lambda|^\gamma$$

for any $\gamma < \min(\sigma+\alpha - 2, \sigma-\frac12, \alpha-\frac12, 1)$.
Equality is possible in the choice of $\gamma$ provided $\sigma,\alpha \not=
\frac32$.
\end{proposition}

\begin{proof} The size conditions in \eqref{eq:Bkernel} guarantees that
$$\norm[K(\lambda)][HS(\sigma,-\alpha)]^2 \les
\lambda^2 \iint_{\{|\x-\y|<\frac1{\lambda}\}} \langle\x\rangle^{-2\sigma}
 \langle\y\rangle^{-2\alpha} d\x d\y \ + \
\iint_{\{|\x-\y|>\frac1{\lambda}\}}\langle\x\rangle^{-2\sigma}\frac1{|\x-\y|^2}
 \langle\y\rangle^{-2\alpha} d\x d\y$$
The first of these integrals is broken up into two domains:
 
Domain 1:  $\max(|\x|,|\y|) \le \frac3{\lambda}$.
 
Domain 2:  $|\y| > \frac2{\lambda},  \x \in B(\y,\frac1{\lambda})$, which also
requires that  $|\x| \sim |\y|$.

\noindent
The second integral is broken up into four domains, namely:
 
Domain 3: $\{|\x-\y| \le \frac12|\y|\}$, which requires $|\y| > \frac2\lambda$
and $|\x| \sim |\y|$.
 
Domain 4:  $\{|\x| \le \frac12|\y|\}$, which requires $|\y| > \frac2{3\lambda}$
 and $|\x-\y| \sim |\y|$.
 
Domain 5:  $\{|\x|, |\x-\y| >\frac12|\y|;\, |\y| > \frac2\lambda\}$, which
requires $|\x| \sim |\x-\y|$.
 
Domain 6:  $\{|\x|, |\x-\y| >\frac12|\y|;\, |\y| \le \frac2\lambda\}$, which
requires $|\x| \sim |\x-\y|$.
In this domain, only values $|\x-\y| > \frac1\lambda > \frac{|\y|}2$ can
make a nonzero contribution.

With the given restrictions on $\sigma$ and $\alpha$ to insure finiteness
of each integral, Domain 1 contributes no more than $C\lambda^{2\gamma}$ to
the total.  Each of the other domains contributes $C\lambda^{2\gamma_i}$,
where $\gamma_i$ is one of the four possible exponents in the definition
of $\gamma$.
\end{proof}

\begin{lemma}
The Fourier transform of $\chi_0(\Bp)'$ in the variable $\lambda$ satisfies
the property
$$\int_{-\infty}^\infty \big\|[\chi_0(\Bp)']^\vee(u)\big\|_{HS(2^+,-2^-)}\, du
< C < \infty$$
uniformly as $\lambda_0 \to 0$.
\end{lemma}
\begin{proof}
First observe that for any pair of points $(\x,\y)$,
$[(\Bp)']^\vee(u)(\x,\y) = \delta(u + |\x - \y|)$, therefore
$$[\chi_0(\Bp)']^\vee(u)(\x,\y) \ = \ \chi_0^\vee(u + |\x-\y|) \ \les \ 
\lambda_0\langle\lambda_0(u+|\x-\y|)\rangle^{-10}$$
The Hilbert-Schmidt norm is bounded above by
\begin{equation*}
\norm[[\chi_0(\Bp)']^\vee(u)][HS(2^+,-2^-)]^2 \les
\lambda_0^2 \iint_{\R^6} \langle\lambda_0(u-|\x-\y|)\rangle^{-20}
\langle\y\rangle^{-4^-}\langle\x\rangle^{-4^-} d\x d\y
\end{equation*}
This is most easily evaluated via the inequality
\begin{equation} \label{eq:shell}
\int_{|\x-\y| = \rho} \langle\x\rangle^{-2\sigma}\,d\x
   \les \langle\rho-|\y|\rangle^{2-2\sigma}
\end{equation}
for $\sigma > 1$.  Integrating with respect to $d\x$ 
over a spherical shell centered at $\y$,
\begin{equation*}
\begin{aligned}
\norm[[\chi_0(\Bp)']^\vee(u)][HS(2^+,-2^-)]^2 \ \les \ 
\lambda_0^2 &\int_0^\infty \int_{\R^3} \langle\lambda_0(u-\rho)\rangle^{-20}
\langle\y\rangle^{-4^-}\langle|\y|-\rho\rangle^{-2^-} d\y d\rho \\
  \les \  &\lambda_0^2\,\int_0^\infty \langle\lambda_0(u-\rho)\rangle^{-20}
\langle\rho\rangle^{-2^-} d\rho
\end{aligned}
\end{equation*}
which leads to the bounds
\begin{equation*}
\norm[[\chi_0(\Bp)']^\vee(u)][HS(2^+,-2^-)] \les \left\{
\begin{aligned}
\lambda_0, \quad &{\rm if}\ u \le \frac2{\lambda_0} \\
\lambda_0^{-9}u^{-10} + \lambda_0^{1/2}u^{-1^-}, \quad &{\rm if}\ 
u \ge \frac2{\lambda_0}
\end{aligned} \right.
\end{equation*}
Integrating this expression yields the quantity 
$C(1 + \lambda_0^{1/2+})$,
which is uniformly bounded as $\lambda_0 \to 0$.
\end{proof}

\begin{lemma} \label{lemma:L1bounds}
The Fourier transform of 
$(\chi_0\Bp)$ in the variable $\lambda$ satisfies the following properties:

\begin{equation} \label{eq:L1bounds}\begin{aligned} 
&\int_{-\infty}^\infty 
\big\|(\chi_0\Bp)^\vee(u)\big\|_{HS(\frac32^+,-1^-)} du  
   < C\lambda_0^{1/2 + }
  \\
&\int_{-\infty}^\infty
\big\|(\chi_0\Bp)^\vee(u)\big\|_{HS(1^+,-\frac32^-)} du 
   < C\lambda_0^{1/2 +}
\end{aligned}
\end{equation}
Identical statements are also true with $\chi_0$ replaced by $\chi_1$.
\end{lemma}

\begin{proof}
First observe that for any pair of points $(\x,\y)$, 
$(\Bp)^\vee(u)(\x,\y) = \frac1{|\x-\y|}\big[\delta(u + |\x - \y|) - \delta(u)
\big]$,
therefore
$$[\chi_0\Bp]^\vee(u)(\x,\y)
    = \frac{\chi_0^\vee (u + |\x-\y|) 
                      - \chi_0^\vee (u)}{|\x-\y|}$$
In the case $|u| \le \frac2{\lambda_0}$,
\begin{equation*}
\Big|[ \chi_0\Bp]^\vee(u)(\x,\y)\Big| \les
  \left\{ \begin{aligned}
  \lambda_0^2&,&\quad &{\rm if}\ |\x-\y| \le \frac1{\lambda_0} \\
  \frac{\lambda_0}{|\x-\y|}&,&\quad &{\rm if}\ |\x-\y| \ge \frac1{\lambda_0}
\end{aligned} \right.
\end{equation*}
where the first estimate comes from the Mean Value theorem.
Up to a factor of $\lambda_0$, this kernel satisfies the hypotheses of 
Proposition \ref{prop:Bnorms}, with the conclusion
$$\norm[[\chi_0\Bp]^\vee(u)][HS(\frac32^+,-1^-)] \les 
  (\lambda_0)^{3/2+}$$
for all $|u| \le \frac2{\lambda_0}$.
In the case $|u| \ge \frac2{\lambda_0}$, we use the fact that $|\chi_0^\vee(u)|
\le \lambda_0\langle\lambda_0 u\rangle^{-10}$ to obtain the pointwise bounds
\begin{equation*}
\Big| [\chi_0\Bp]^\vee(u)(\x,\y) \Big| \les
   \left\{ \begin{aligned}
 &\frac{1}{\lambda_0^8|u|^{10}},&\quad &{\rm if}\ |\x-\y| \le \frac{|u|}2 \\
  &\frac1{\lambda_0^8|u|^9|\x-\y|},&\quad &{\rm if}\ |\x-\y| \ge 2|u| \\
 \frac{\lambda_0}{|u|}\langle\lambda_0(&u + |\x-\y|)\rangle^{-10},
&\quad &{\rm if}\ \frac{|u|}2 < |\x-\y| < 2|u|
 \end{aligned} \right.
\end{equation*}
The restriction of this kernel to the domain $\{|\x-\y| \ge 2|u|\} \cup
\{|\x-\y| \le \frac{|u|}2\}$ must have
$HS(\frac32^+,-1^-)$-norm of $\lambda_0^{-8}|u|^{-9.5^-}$, also by Proposition
\ref{prop:Bnorms}.  Since we are assuming $|u| \ge \frac2{\lambda_0}$,
this is less than $|u|^{-1.5^-}$.
To estimate the Hilbert-Schmidt
norm of the remaining annular piece, we once again use the inequality
\begin{equation} \tag{\ref{eq:shell}}
\int_{|\x-\y| = \rho} \langle\x\rangle^{-2\sigma}\,d\x 
   \les \langle\rho-|\y|\rangle^{2-2\sigma}
\end{equation}
for $\sigma > 1$.  Thus
\begin{multline*}
\frac{\lambda_0^2}{|u|^2}\iint_{|\x-\y|\sim |u|} \langle\x\rangle^{-3^-}
  \langle\lambda_0(u + |\x-\y|)\rangle^{-20} \langle\y\rangle^{-2^-}d\x d\y
\\
\les \frac{\lambda_0^2}{|u|^2}\int_{\frac{|u|}2}^{2|u|} 
    \langle\lambda_0(u + \rho)\rangle^{-20} \int_{\R^3} 
    \langle\rho-|\y|\rangle^{-1^-}\langle\y\rangle^{-2^-}d\y d\rho
\\
\les \frac{\lambda_0^2}{|u|^2} \int_{\frac{|u|}2}^{2|u|}  
   \langle\rho\rangle^{-\epsilon} \langle\lambda_0(u+\rho)\rangle^{-20} d\rho
\les \frac{\lambda_0}{|u|^{2+\epsilon}}
\end{multline*}
Putting the pieces together, it follows that
$$\norm[[\chi_0\Bp]^\vee][HS(\frac32^+,-1^-)]
   \les (\lambda_0)^{3/2+\epsilon} \langle\lambda_0 |u|\rangle^{-1-\epsilon}$$
proving the first claim of the lemma.  The second line of equation 
\eqref{eq:L1bounds} follows from symmetry in the variables $\x$ and $\y$.
\end{proof}

Recall from equation \eqref{eq:Bprimeseries} that 
\begin{equation*}
\chi_0(\lambda)(\Btp)'(\lambda) \ = \ \sum_{n=1}^\infty\sum_{m=0}^{n-1} (-1)^n
\big((\chi_1\Bp)(\lambda)VS_0^{-1}\big)^m
\chi_0(\lambda)(\Bp)'(\lambda)VS_0^{-1}
\big((\chi_1\Bp)(\lambda)VS_0^{-1}\big)^{n-(m+1)}
\end{equation*}

Multiplication of operator-valued functions results in a convolution of their
respective Fourier transforms, just as it does in the scalar case.  Similarly,
the $L^1$ theory of convolution applies in this setting provided the domain
of each operator is identified with the range of its predecessor.  
Then
\begin{multline*}
\norm[[\chi_0(\Btp)']^\vee][L^1(HS(-1^-,-2^-))] \\ 
\le  \sum_{n=1}^\infty \sum_{m=0}^{n-1}
\norm[(\chi_1\Bp)^\vee VS_0^{-1}][]^m
\norm[[\chi_0(\Bp)']^\vee VS_0^{-1}][]
\norm[(\chi_1\Bp)^\vee VS_0^{-1}][]^{n-(m+1)} 
\end{multline*}    
where the norms are taken in $L^1$ with values in $HS(-2^-,-2^-)$, 
$HS(-1^-,-2^-)$, and $HS(-1^-,-1^-)$, respectively.

From equation \eqref{eq:L1bounds} we see that the sum converges exponentially
provided $\lambda_0$ is chosen small enough.

\smallskip \noindent
{\bf Remarks.} Throughout the discussion, operators have been estimated by the 
Hilbert-Schmidt norm as a matter of computational convenience. More precisely,
we needed to know that various kernels $K(x,y)$ of $L^2$-bounded operators 
have the property that $|K(x,y)|$ again gives rise to an $L^2$-bounded operator
(on this level of generality we do not need to distinguish between $L^{2,\sigma}$
and $L^2$, since the weights can be included in the kernel).
Note that this property is automatic if $K(x,y)$ is Hilbert-Schmidt.

More generally, note that an operator of the form $I+T$ with kernel
$\delta(x-y)+K(x,y)$ where $K$ is Hilbert-Schmidt, still has the property
that the absolute value of the kernel gives rise to an $L^2$-bounded operator.
Moreover, if $T$ is Hilbert-Schmidt and $(I+T)^{-1}$ exists, then $(I+T)^{-1}-I=
-(I+T)^{-1}T$ is again Hilbert-Schmidt. This observation implies, in particular,
that $|S_0^{-1}|$ is $L^{2,\sigma}$-bounded with $-\frac52<\sigma<-\frac12$. 
Here and in what follows, $|T|$ stands for the 
operator that is given by the absolute value of the kernel of $T$.

Consider the case when the derivative $\frac{d}{d\lambda}$ falls on a
different term in \eqref{eq:IBP}, for example on $e^{i\lambda|\y-x_4|}$.
In the lines which follow, one is then led to control the size of
\begin{equation*}
 \sup_{\x,\y\in\R^3} \frac1{t^{3/2}}
\norm[V][L^{2,1^+}]
 \ \big\|{\textstyle \int} |S_0^{-1} [\chi_0(\Btp)]
  ^\vee(u)|du \big\|_{{\mathcal B}(-1^-,-1^-)} \ \big\||\x-\cdot|^{-1}\big\|_{L^{2,-1^-}}
\end{equation*}
which depends eventually on the finiteness of the central integral
$$\int_{-\infty}^\infty \norm[|{[\chi_0 \Btp]}^\vee(u)|][{\mathcal B}(-1^-,-1^-)]\, du.$$
Here ${\mathcal B}(-1^-,-1^-)$ stands for the bounded operators $L^{2,-1-}\to L^{2,-1-}$. 

Unlike in \eqref{eq:Bprimeseries}, the Neumann series for $\chi_0\Btp$ 
begins with a zero-order term, namely $\chi_0(\lambda)$ times the identity map.
While the identity is a bounded operator on $L^{2,-1^-}$ it does not
belong to the Hilbert-Schmidt class.  All higher order terms are 
Hilbert-Schmidt, however, because they each contain at
least one multiple of $\Bp(\lambda)VS_0^{-1}$.
A similar zero-order term appears if the derivative in \eqref{eq:IBP}
falls on $e^{i\lambda|\x-x_1|}$ or on the cut-off function $\chi_0(\lambda)$.


{\bf Acknowledgement:} The second author was supported by the NSF grant
DMS-0070538 and a Sloan fellowship.

\bibliographystyle{amsplain}

\medskip\noindent
\textsc{Division of Astronomy, Mathematics, and Physics, 253-37 Caltech, Pasadena, CA 91125, U.S.A.}\\
{\em email: }\textsf{\bf mikeg@caltech.edu, schlag@caltech.edu} 

\end{document}